\documentclass[12pt]{amsart}
\usepackage{amsmath, amssymb, amsthm, latexsym}
\usepackage{slashed}
\input xypic
\newtheorem{theorem}{Theorem}
\newtheorem{proposition}[theorem]{Proposition}

\newtheorem{lemma}[theorem]{Lemma}

\newtheorem{corollary}[theorem]{Corollary}
\newtheorem{definition}[theorem]{Definition}

\def\Proof{\medskip\noindent{\bf Proof: }}

\def\Z{\mathbb{Z}}

\def\C{\mathbb{C}}

\def\R{\mathbb{R}}
\def\C{\mathbb{C}}

\def\N{\mathbb{N}}

\def\Pi{\mathbb{P}^{\infty}}

\def\mb{\mathfrak{B}}

\def\qed{\hfill$\square$\medskip}

\def\Zpk{\mathbb{Z}/p^{k}}
\def\Zpk1{\mathbb{Z}/p^{k-1}}

\newcommand{\rref}[1]{(\ref{#1})}

\newcommand{\cform}[3]{\begin{array}{c}
{\scriptstyle #3}\\
#1\\
{\scriptstyle #2}\end{array}}

\newcommand{\beg}[2]{\begin{equation}\label{#1}#2\end{equation}}
\def\r{\rightarrow}

\def\mb{\mathcal{B}}

\def\sl2{\widetilde{SL_{2}(\Z)}}

\begin{document}

\title{Categorical geometry and integration without points}
\author{Igor Kriz and Ale\v{s} Pultr}
\thanks{The first author was supported in part by NSF grant DMS 1104348.
The second author was supported by the projects 
1M0545 and MSM 0021620838 of the Ministry of Educations of the 
Czech Republic}
\subjclass[2010]{28A60,16B50,03G30,28C20}
\keywords{point free measures, boolean rings, categorical geometry,
locales, Segal space}
\dedicatory{{\em(}dedicated to the memory of Irving Segal{\em)}}
\begin{abstract}
The theory of integration over infinite-dimensional spaces is known to encounter 
serious difficulties. Categorical ideas seem to arise naturally on
the path to a remedy. 
Such an approach was suggested and initiated by Segal in his pioneering article \cite{segal}.
In our paper we follow his ideas from a different perspective, slightly more categorical, and strongly inspired by the point-free topology.

First, we develop a general (point-free) concept of measurability (extending the standard Lebesgue integration when applying to the classical $\sigma$-algebra). Second (and here we have a major difference from the classical theory), we prove that every
finite-additive function $\mu$ with values in $[0,1]$ can be extended to a measure on an abstract $\sigma$-algebra; this correspondence is functorial and yields uniqueness. As an example we show that the Segal space can be characterized by completely canonical data. Furthermore, from our results it follows that a satisfactory point-free integration arises everywhere where we have  a finite-additive probability function on a Boolean algebra.
\end{abstract}
\maketitle

\section{Introduction}

\vspace{3mm}
The basic concept of a $\sigma$-algebra, meaning a system of subsets of
a given set closed under complements and countable unions,
and the accompanying concept of $\sigma$-additive measure \cite{leb},
are the cornerstones of the modern theory of integration. 
Yet, these concepts are tested to the extreme (and sometimes beyond) in contexts 
of stochastic analysis \cite{mal} on the one hand, and quantum field theory 
(see \cite{ias} for a relatively recent introduction) on the
other. In these cases, one often desires a theory of integration over 
infinite-dimensional spaces where $\sigma$-additive measures are either not
present or difficult to construct. A classical example
is provided by the case of attempting to define a Gaussian measure on 
the $\sigma$-algebra generated by cylindrical sets \cite{kuo, mal}
in an infinite-dimensional Hilbert space $H$.
This is impossible, because the ``Gaussian measure'' on cylindrical sets themselves
is not $\sigma$-additive. One known solution of this problem \cite{radon,radon1}
uses a method of 
``enlarging'' the space $H$. This became known as ``Radonification''. Its most
famous example is the Wiener measure \cite{wiener, wiener1}. 
The disadvantage of this method is that the ``enlargement'' is intrinsically
non-canonical (in particular not functorial in the Hilbert space with respect
to isometries). A different solution \cite{grosssegal} became known as the
so called Segal space. This construction can be made functorial, but
still involves an arbitrary $\sigma$-algebra: functoriality is achieved
by forgetting a part of the data. 

\vspace{3mm}

Irving Segal, actually, very much believed in the idea of {\em algebraic integration}, 
i.e. a theory of integration which
would start with an algebraic structure {\em without} an underlying 
space of points. In his 1965 discoursive article \cite{segal}, he
laid out the perspectives of such a theory around that time. 
At the beginning of the article, he mentions early efforts of
describing Lebesgue integration 
in terms of Boolean rings by Carath\'{e}odory
and his associates, and the ultimate snag of this direction due
to the fact that functions only enter that picture in a circumlocutory 
fashion, forestalling a comprehensive development of integration
theory along these lines.
He went on to describe many 
more sophisticated efforts to define point-free integration in more
modern times. A prominent example is an approach to quantum field theory
of completing the algebra
of bounded continuous cylinder functions on a Hilbert space $H$,
completing it with respect to the $sup$-norm, and applying the Gelfand
construction and the Riesz representation theorem to produce a compact Hausdorff
space $X$. The 1987 book by Glimm and Jaffe \cite{gj} describes many applications
and developments in that direction. 
In \cite{segal}, Segal also goes on to describe other
``point-free'' developments in integration theory, including his own 1954 paper
\cite{segalkol} on the Kolmogoroff theorem, and integration theories based
on Von Neumann algebras. Such efforts continued to more modern
times: to give just one example, in the 2002 paper \cite{coquand}, T.Coquand and E.Palmgren
described a Daniell-type integration, considering
measures on Boolean algebras with a ``strong apartness 
relation'', and extending them to a metric completion.

\vspace{3mm}
The space $X$, however, is hard to get
one's hands on in the context of analysis. The
large mathematical
physics community preferred to use a non-canonical Schwartz distribution space to support
a countably additive version of this measure rather than that compact Hausdorff space. 
Leonard Gross invented the notion of Abstract Wiener space in order to develop 
potential theory on a Hilbert space $H$. The norm on the
abstract Wiener space, invented by L.Gross, is 
essential for formulating some 
natural theorems in potential theory. Thus, we come full circle, and points seem
to force their way back into the theory.

\vspace{3mm}
In this paper, we approach Segal's idea of `point-free integration'
from a different, surprisingly naive and fundamental point of view, going,
in a way, back to Carath\'{e}odory, but with a new take
provided by category theory. During the first author's Prague years, the 
authors of this paper collaborated on a series of papers in category theory,
developing analogues of certain constructions and
facts of point set topology in the point-free context.
In the theory of locales, a topological space is replaced by 
a ``completely $\bigvee-\wedge$ distributive lattice'' (frame)
which models the algebraic structure of the set of 
open sets. Point-free topology has a long history. 
It boomed since the late fifties 
(out of the many authors and articles from 
the early time let us mention e.g.
\cite{Ehr} and \cite{IsbAPS}). The development of 
the first decades culminated in Johnstone's monograph
\cite{john} (see also his excellent surveys \cite{Jpoint} 
and \cite{Jart}), and continues still, in the more recent 
decades particularly in the theory of enriched point-free 
structures. Both authors also published in the field 
(e.g. \cite{kriz, kriz1, pultr}). For more about frames 
and further references see e.g. \cite{john}, \cite{Jart},
\cite{Jpoint}, \cite{pp}, \cite{pultr}.
It should be mentioned that our investigations in 
this paper are also related to a special branch 
of pointfree topology, the theory of $\sigma$-frames 
(motvated by the so called Alexandroff spaces) where one
assumes countable suprema and ``countable $\bigvee-\wedge$ 
distributivity'' only (see e.g. \cite{Gil} or \cite{BBGil}).

\vspace{3mm}
The lesson of point-free topology is that techniques of category theory can often
be used to supply concepts which seemingly need points.
In the context of this paper, the key point is that we
find a categorical version of the appropriate concept of a {\em measurable
function}, since such objects can be no longer defined in the present context
by values on points. Our definition generalizes
the notion of measurable function in the case of a classical $\sigma$-algebra,
and a theory of integration which gives the same results as Lebesgue integration
in this classical case.
The first main result of the present paper is that measure and Lebesgue integration
theory can indeed be generalized to the context of abstract $\sigma$-algebras.

\vspace{3mm}
The second main result 
of this paper (Theorem \ref{trad}) shows a major
difference from the classical
theory: For every finite-additive function $\mu$ on a Boolean algebra $\mathcal{B}$ with
values in $[0,1]$,
there exists a measure on an abstract $\sigma$-algebra extending $\mu$. Moreover,
by using appropriate conditions,
one can obtain {\em uniqueness} (and hence functoriality).
We show as an example how from a point-free point of view, the Segal space can be
characterized by data completely canonical (up to canonical isomorphism). More generally,
by the results of the present paper, a point-free measure, 
and point-free integration theory arises everywhere where we have a finite-additive
probability function on a Boolean algebra.

\vspace{3mm}
To place the present paper in the context of existing work,
we refer to the ultimate authority, namely D.H. Fremlin's book
\cite{dhf}. While that book documents many ideas relevant to
the topic of the present paper, astonishingly, the basic
and naive geometric view which we present here seems to have
been bypassed by the field. The Stone representation theorem
\cite{dhf}, p.70, asserts that every abstract $\sigma$-algebra
with a measure is isomorphic to the $\sigma$-algebra of some
measure space. Maharam's theorem (Chapter 33 of \cite{dhf})
gives yet another characterization of the possible 
isomorphism types of abstract $\sigma$-algebras with measure
(under some mild assumptions). Although this is not
stated explicitly in \cite{dhf}, the statement of our Theorem
\ref{trad} can actually be deduced from the material
covered there (see \cite{dhf}, Exercise
325Y (b) on p. 103). All these results use complicated
pointed measure constructions as intermediaries. 
In contrast, our direct categorical ``Grothendieck-style''
treatment of abstract $\sigma$-algebras develops 
point-free integration as an intrinsic geometry (in analogy,
for example, with the intrinsic categorically-geometric
treatment of super-manifolds in 
\cite{ias} Chapter 3, which is of great 
benefit even though
all super-manifolds are in fact ``spatial''). 

\vspace{3mm}
The present paper is organized as follows: the development of integration theory
for measures on abstract $\sigma$-algebras is done in Section \ref{smeasure}.
The extension of finite-additive to $\sigma$-additive measures is
handled in Section \ref{sradon} below.

\vspace{3mm}
{\bf Acknowledgement:} The authors are very indebted to Adelchi Azzalini
and Leonard Gross for valuable discussions.

\section{Measure and integration on abstract $\sigma$-algebras}

\label{smeasure}

\begin{definition}
\label{dcb}
By an  {\em abstract $\sigma$-algebra}, we mean a Boolean algebra $B$
in which there exist countable joins.
By a morphism of abstract $\sigma$-algebras we mean a map which preserves order,
$0,1$, complements, and countable joins (and therefore meets).
\end{definition}

\vspace{3mm}
\noindent
This is an obvious definition, but let us review a few useful facts.
We automatically have
the distributivity
\beg{edcb}{\left(\cform{\bigvee}{i=0}{\infty}a_i\right)\wedge b=
\cform{\bigvee}{i=0}{\infty}(a_i\wedge b).
}
This is because the operation $(\smallsetminus a)\vee ?$ is right
adjoint to $a\wedge ?$ considered as self-functors of the POSET
$B$. On the other hand, it is {\em not} automatic for morphisms
of Boolean algebras $B\r C$ where $B,C$ are abstract $\sigma$-algebras to
preserve countable joins (take, for example, the morphism from
the Boolean algebra of all subsets on $\N$ to $\{0,1\}$ given
by an ultrafilter which does not correspond to a point.

Recall that Boolean algebras can be characterized as commutative
associative unital rings satisfying the relation
$$a^2=a.$$
In this identification, the product corresponds to the meet, and $+$ to
symmetric difference $\oplus$. (Note that since join and meet
are symmetric in a Boolean algebra, another symmetric identification
is possible; however, this is the usual convention.)
From the point of view of this identification, the coproduct of
Boolean algebras is the tensor product. 

Abstract $\sigma$-algebras can be similarly characterized as universal algebras with operations of
at most countable arity, for example as commutative associative
unital rings 
satisfying the identity 
$$a^2=a$$
with an operation of countable arity 
$$\cform{\prod}{i=0}{\infty} a_i$$
which satisfies infinite commutativity, infinite associativity
$$\cform{\prod}{i=0}{\infty} \left( \cform{\prod}{j=0}{\infty} a_{ij}\right)
=\cform{\prod}{i,j=0}{\infty} a_{ij},$$
$$\cform{\prod}{i=0}{\infty} a_i=\cform{\prod}{i=0}{n} a_i\;\text{when $a_{n+1}=
a_{n+2}=...=1$},$$
and 
$$\cform{\prod}{i=0}{\infty} (a_i+a_ib+b)=\cform{\prod}{i=0}{\infty} a_i +
b\cform{\prod}{i=0}{\infty}a_i
+b.$$
From this point of view, the coproduct of 
abstract $\sigma$-algebras $B_1$, $B_2$
can also be characterized
as a kind of ``completed tensor product'', concretely the quotient of the free 
abstract $\sigma$-algebra on
the set
$$\{b_1\otimes b_2\;|\; b_i\in B_i\}$$
by the relations
$$\cform{\bigvee}{n=0}{\infty}(a_n\otimes b)\sim \left(
\cform{\bigvee}{n=0}{\infty} a_n\right)\otimes b,$$
$$\cform{\bigvee}{n=0}{\infty}(b\otimes a_n)\sim \left(
b\otimes\cform{\bigvee}{n=0}{\infty} a_n\right),$$
$$0\otimes a\sim a\otimes 0\sim 0.$$

\vspace{3mm}
In this paper, we will describe a theory of measure and Lebesgue integration
on abstract $\sigma$-algebras. A {\em measure} on an 
abstract $\sigma$-algebra (resp. an {\em additive
function} on a Boolean algebra) $B$ is a function
$$\mu:B\r [0,\infty]$$
which satisfies $\mu(0)=0$ (in both cases) and
$$\mu\left(\cform{\bigvee}{i=0}{\infty} a_i\right)=
\cform{\sum}{i=0}{\infty} \mu(a_i)$$
when $i\neq j\Rightarrow a_i\wedge a_j=0$ (resp.
$$\mu(a\vee b)=\mu(a)+\mu(b)$$
when $a\wedge b=0$).

An abstract $\sigma$-algebra 
(resp. Boolean algebra) $B$ with a measure (resp. additive
function $\mu$) is called {\em reduced} if
$$\mu(a)=0\Rightarrow a=0.$$
An {\em ideal} $I$ on an 
abstract $\sigma$-algebra (resp. Boolean algebra $B$) is 
a subset $I\subseteq B$ such that $0\in I$,
when $a\in I$ and $b\leq a$ then $b \in I$ (in both cases) and
$$\cform{\bigvee}{i=0}{\infty} a_i\in I \;\text{if $a_i\in I$}$$
(resp. 
$$a\vee b\in I\;\text{if $a,b\in I$}.$$

\begin{lemma}
Let $B$ be an 
abstract $\sigma$-algebra (resp. a Boolean algebra). Let $I\subseteq B$
be an ideal. Define an equivalence relation $\sim$ on $B$ by
$$a\sim b \;\text{if $a\oplus b\in I$.}$$
Then $B/I:=B/\!\!\sim$ with operations induced from $B$ is an 
abstract $\sigma$-algebra
(resp. Boolean algebra).
\end{lemma}
\qed

\vspace{3mm}
\noindent
{\bf Example:} 
If $\mu$ is a measure (resp. additive function) on an 
abstract $\sigma$-algebra
(resp. Boolean algebra) $B$, then
$$I_\mu:=\{a\in B\;|\;\mu(a)=0\}$$
is an ideal and $B/I_\mu$ has an induced measure (resp. additive
function) from $B$ with respect to which it is reduced.

\vspace{3mm}
By a {\em $\sigma$-frame} we shall mean a lattice which contains $0,1$, finite
meets and countable joins which are distributive. By a {\em universal 
abstract $\sigma$-algebra}
on a $\sigma$-frame we shall mean the left adjoint of the forgetful functor
from abstract $\sigma$-algebras to $\sigma$-frames.

\begin{lemma}
\label{lb}
The following abstract $\sigma$-algebras are canonically isomorphic:
\begin{enumerate}

\item
\label{i1}
The quotient of the free 
abstract $\sigma$-algebra on the set of intervals $[0,t]$, 
$0\leq t\leq\infty$ by the relations $[0,s]\leq [0,t]$ for $s\leq t$ and
$\bigwedge [0,s_n]= [0,t]$ when $s_n\searrow t$

\item
\label{i2}
The universal abstract $\sigma$-algebra on the $\sigma$-frame of open sets
with respect to the analytic topology
on $[0,\infty]$.

\item
\label{i3} 
The quotient of the free 
abstract $\sigma$-algebra on the set of intervals $[t,\infty]$, 
$0\leq t\leq\infty$ by the relations $[t,\infty]\leq [s,\infty]$ for $s\leq t$ and
$\bigwedge [s_n,\infty]= [t,\infty]$ when $s_n\nearrow t$

\end{enumerate}

\end{lemma}

\Proof
Using universality, one constructs maps of 
abstract $\sigma$-algebras between 
\rref{i1} and \rref{i2} in both directions:
a map from \rref{i1} to \rref{i2} is given by
$$[0,t]\mapsto \smallsetminus (t,\infty].$$
A map from \rref{i2} to \rref{i1} is given by
$$(s,t)\mapsto (\bigvee [0,t_i])\smallsetminus [0,s]$$
where $t_i\nearrow t$.
Correctness follows from compactness of closed intervals in $[0,\infty]$.
The maps are obviously inverse to each other. Inverse maps of
abstract $\sigma$-algebras 
between \rref{i2} and \rref{i3} are constructed analogously.
\qed

\vspace{3mm}
Let $\mb$ be the 
abstract $\sigma$-algebra characterized by the equivalent properties of
Lemma \ref{lb}. 

\vspace{3mm}
\noindent
{\bf Problem:} Is $\mb$ isomorphic to the abstract $\sigma$-algebra of
Borel sets in $[0,\infty]$?

\vspace{3mm}
We will also denote by $\mb(\R)$ (resp. $\mb(\C)$) the universal 
abstract $\sigma$-algebras on the $\sigma$-frame of 
open sets in $\R$ (resp. $\C$) with
respect to the analytic topology.
A {\em non-negative measurable function} on an 
abstract $\sigma$-algebra $\Sigma$
is a morphism of abstract $\sigma$-algebras
$$f:\mb\r\Sigma.$$
For measurable functions $f,g:\mb\r \Sigma$, we write
$$f\leq g$$
if for every $t\in[0,\infty]$, 
$$f[0,t]\geq g[0,t].$$
Addition
$$+:[0,\infty]\times [0,\infty]\r [0,\infty]$$
induces a map of abstract $\sigma$-algebras 
$$+:\mb\r\mb\amalg\mb$$
(here $\amalg$ denotes categorical coproduct).
The {\em sum} of measurable functions $f,g:\mb\r\Sigma$
is the composition
$$\diagram
\mb\rto^+&\mb\amalg\mb\rto^{f\amalg g}&\Sigma\amalg\Sigma\rto^{\nabla} & \Sigma
\enddiagram$$
where $\nabla$ is the categorical co-diagonal. A similar construction clearly
works with $+$ replaced by any continuous operation on $[0,\infty]$, such as
multiplication, or the $\min$ or $\max$ function.

Let $S$ be a set. Denote by $2^S$
the abstract $\sigma$-algebra of all subsets of $S$.

\begin{lemma}
\label{l2s}
Every morphism of abstract $\sigma$-algebras $\phi:B\r 2^{\{0\}}$ is induced
by a point in $[0,\infty]$, i.e. a map of sets $\{0\}\r[0,\infty]$.
\end{lemma}

\Proof
Clearly, we cannot have $\phi(\{x\})=\phi(\{y\})=1$ for
two points $x\neq y$, since then
$$0=\phi(0)=\phi(\{x\}\wedge\{y\})=\phi(\{x\})\wedge\phi(\{y\})=1.$$
Thus, there is at most one point $x\in[0,\infty]$ with 
\beg{ei1}{\phi(\{x\})=1.}
We claim that an $x$ satisfying \rref{ei1} exists. In effect, assuming
$\phi(\{\infty\})=0$, we must have $\phi[n_i/2^i,(n_i+1)/2^i)=1$
for some $n_i\in \N$ by countable additivity. Then 
$$\phi\left(\cform{\bigwedge}{i=0}{\infty} [n_i/2^i,(n_i+1)/2^i)\right)=1,$$
but the argument on the left hand side has at most one point. On the other
hand, it cannot be
empty because $\phi(0)=0$.
\qed

\vspace{3mm}
A {\em simple function} on an 
abstract $\sigma$-algebra $\Sigma$ is a
morphism of abstract $\sigma$-algberas $F:\mb\r \Sigma$ for which there exists
a finite Boolean algebra (hence abstract $\Sigma$-algebra)  $\mathcal{F}=2^{\{0,...,n\}}$
and morphisms of abstract $\sigma$-algberas $\chi:\mb\r\mathcal{F}$, $s:\mathcal{F}\r \Sigma$
such that the diagram
$$\diagram
\mb\drto_f\rto_\chi &\mathcal{F}\dto^s\\
&\Sigma
\enddiagram
$$
commutes.

By Lemma \ref{l2s}, $\chi$ is induced by
$$i\mapsto x_i\in [0,\infty],\; i=0,...,n.$$
Let $\mu:\Sigma\r[0,\infty]$ be a measure. Then define
$$\int fd\mu:=\cform{\sum}{i=0}{n}\mu(s(\{i\})x_i.$$
For an arbitrary function $f:\mb\r\Sigma$, define
$$\int fd\mu=sup\left\{\int gd\mu\;|\;g\leq f \;\text{and $g$ is simple}\right\}.$$

\vspace{3mm}
\noindent
{\bf Remark:}
Clearly, when $(X,\Sigma,\mu)$ is a measurable space, $\phi:X\r[0,\infty]$ and
is (Borel) measurable, then defining $f:\mb\r\Sigma$ by
$$f(a):=\phi^{-1}(a),$$
we have by definition
$$\int fd\mu=\int\phi d\mu.$$

\vspace{3mm}
\begin{lemma}
\label{ls1}
Let $\Sigma$ be an 
abstract $\sigma$-algebra with measure $\mu$. If $s_1,s_2:\mb\r\Sigma$
are simple functions, then $s_1+s_2$ is a simple function and
$$\int(s_1+s_2)d\mu=\left(\int s_1d\mu\right)+\left(\int s_2d\mu\right).$$
\end{lemma}
\qed

\vspace{3mm}
For an 
abstract $\sigma$-algebra 
$\Sigma$ and functions $f,f_n;\mb\r \Sigma$, $n\in\N$,
we write $f_n\nearrow f$ if 
$$f_0\leq f_1\leq f_2\leq...\leq f$$ 
and for all $t\in [0,\infty]$,
$$f([0,t])=\cform{\bigwedge}{n}{}f_n([0,t]).$$
(Note: The non-trivial inequality is $\leq$.) It is worth noting that
by characterization \rref{i1} of Lemma \ref{lb},
whenever
$$f_0\leq f_1\leq f_2\leq...$$
for non-negative measurable functions $f_i:\mb\r\Sigma$, 
there exists a unique non-negative measurable function $f:\mb\r \Sigma$ such that
$f_i\nearrow f$. Analogously, one defines for non-negative measurable functions
$$f_0\geq f_1\geq f_2....$$
$f_n\searrow f$ when 
$$f([t,\infty])=\cform{\bigwedge}{n}{}f_n([t,\infty]).$$
Using characterization \rref{i3} of Lemma \ref{lb},
one shows that such non-negative measurable function $f$ always exists.

\vspace{3mm}
\begin{lemma}
\label{ls2}
For an 
abstract $\sigma$-algebra $\Sigma$ and a 
non-negative measurable function $f:\mb\r\Sigma$,
there exist simple functions $s_n:\mb\r\Sigma$ such that $s_n\nearrow f$.
\end{lemma}

\Proof
Let
$$S_n(k):=\left[\frac{k}{2^n},\frac{k+1}{2^n}\right)\;\text{if $k=0,...,(n2^n)-1$},$$
$$S_n(n2^n):=[n,\infty].$$
Then let, for $a\in\mb$,
$$s_n(a):=\cform{\bigvee}{\{k\in\{0,...,n2^n\}|k/2^n\in a\}}{}f(S_n(k)).$$
\qed

\vspace{3mm}
\begin{theorem}
\label{tm}
(Levi's theorem, Lebesgue monotone convergence theorem) let $\Sigma$ be an 
abstract $\sigma$-algebra
with measure $\mu$ and let $f, f_n:\mb\r \Sigma$ be functions, 
$f_n\nearrow f$. Then
$$\cform{\lim}{n\r \infty}{}\int f_nd\mu=\int fd\mu.$$
\end{theorem}

\Proof
The $\leq$ inequality is immediate. Now for a simple function $s\leq f$, let
$$\diagram\mb\rto^g\drto_s &2^{\{0,...,n\}}\dto^{\tilde{s}}\\
&\Sigma\enddiagram$$
let $\tilde{s}(i):=a_i$ and let $x_i\in [0,\infty]$ be such that
$$g(a)=\{i\;|\;x_i\in a\}$$
(see Lemma \ref{l2s}). Let 
$$b_{ink}:=f_n([y_{ik},\infty])$$
where $y_{ik}=x_i-1/k$ if $x_i\neq\infty$, and $y_{ik}=k$ if $x_i=\infty$.
Then $f_n\nearrow f$ implies
\beg{etm+}{\cform{\bigvee}{n}{}b_{ink}\leq a_i.
}
This implies
$$\cform{\lim}{n\r\infty}{}\int f_nd\mu
\geq \cform{\lim}{n\r\infty}{}\cform{\sum}{i=0}{m}\mu(b_{ink}\wedge a_i)y_{ik}
=\cform{\sum}{i=0}{m}\mu(a_i)y_{ik}.$$
By \rref{etm+}, the limit of the right hand side as $k\r \infty$ is 
$$\int sd\mu.$$
Since $s\leq f$ was an arbitrary simple function, the $\geq$ inequality follows.
\qed

\vspace{3mm}

\begin{lemma}
\label{l+1}
If, for an 
abstract $\sigma$-algebra $\Sigma$ and functions $f_1,f_2,g_1,g_2:\mb\r\Sigma$,
$f_1\leq f_2$ and $g_1\leq g_2$, then
$$f_1+f_2\leq g_1+g_2.$$
\end{lemma}

\Proof
Follows from the definition and the fact that there exists a countable
sequence of pairs $0\leq r_i,s_i$, $r_i+s_i\leq t$, such that
$$+([0,t])=\cform{\bigvee}{i}{}[0,r_i]\otimes [0,s_i].$$
\qed

\vspace{3mm}

\begin{lemma}
\label{l+2}
For an
abstract $\sigma$-algebra $\Sigma$ and non-negative measurable functions
$f,f_n:\mb\r\Sigma$, $n\in \N$, if we have $f_n\nearrow f$,
$g_n\nearrow g$, then $f_n+g_n\nearrow f+g$.
\end{lemma}

\Proof
$$f_0+g_0\leq f_1+g_1\leq...\leq f+g$$
follows from Lemma \ref{l+1}. But there exist 
$0\leq r_{i1},...,r_{in_i},s_{i1},...,s_{in_i}$ such that
\beg{el+2+}{+([0,t])=\cform{\bigwedge}{i=0}{\infty}
\cform{\bigvee}{j=1}{n_i}[0,r_{ij}]\otimes[0,s_{ij}].
}
Then
\beg{el+21}{(f\amalg g)(S)=\cform{\bigwedge}{n}{}(f_n\amalg g_n)(S)
}
is true for $S=[0,r_{ij}]\otimes [0,s_{ij}]$, hence by distributivity
for
$$S=\cform{\bigvee}{j=1}{n_i}[0,r_{ij}]\otimes [0,s_{ij}].$$
To see this, compute
$$(f\amalg g)(S)=
\cform{\bigvee}{j=1}{n_i}\cform{\bigwedge}{n}{}(f_n\amalg g_n)
([0,r_{ij}]\otimes [0,s_{ij}])$$
$$=\cform{\bigwedge}{m_1,...,m_{n_i}}{}
\cform{\bigvee}{j=1}{n_i}(f_{m_j}\amalg g_{m_j})([0,r_{ij}]\otimes [0,s_{ij}])$$
$$=\cform{\bigwedge}{n}{}(f_n\amalg g_n))(S).$$
Hence, \rref{el+21} holds for $S=+([0,t])$ by \rref{el+2+}.
\qed

\vspace{3mm}
\begin{theorem}
\label{tadd}
For an
abstract $\sigma$-algebra 
$\Sigma$ with measure $\mu$ and $f,g:\mb\r\Sigma$ non-negative
measurable functions,
$$\int(f+g)d\mu=\int fd\mu +\int gd\mu.$$
\end{theorem}

\Proof
By Lemma \ref{ls2}, choose simple functions $s_n\nearrow f$, $s^{\prime}_{n}\nearrow g$.
We have
$$\int(f+g)d\mu = \;\text{by Lemma \ref{l+2} and Lebesgue's theorem \ref{tm}}$$
$$=\cform{\lim}{n\r\infty}{}\int(s_n+s^{\prime}_{n})d\mu = \;\text{by Lemma \ref{ls1}}$$
$$=\cform{\lim}{n\r\infty}{}\int s_nd\mu + \cform{\lim}{n\r\infty}{}\int s^{\prime}_{n}d\mu=
\;\text{by theorem \ref{tm}}=$$
$$=\int fd\mu+\int gd\mu.$$
\qed

\vspace{3mm}
Since we showed that minima, maxima, and limits of increasing 
and decreasing sequences of non-negative measurable functions are defined
as non-negative measurable functions, one can define $\liminf$ and $\limsup$
among non-negative measurable funcions on an abstract $\sigma$-algebra.

\vspace{3mm}
\begin{lemma}
(Fatou's lemma)
If $f_n:\mb\r \Sigma$ are non-negative measurable functions for $n\in\N$,
then
\beg{ef1}{\int\left(\cform{\liminf}{n\r\infty}{}f_n\right)d\mu
\leq \cform{\liminf}{n\r\infty}{}\int f_nd\mu.
}
\end{lemma}

\Proof
(following \cite{rudin}): Put
$$g_k:=\cform{\inf}{i\geq k}{}f_i,$$
and apply Lebesgue monotone convergence theorem to $g_k$.
\qed

\vspace{3mm}

Now for an abstract $\sigma$-algebra $\Sigma$, by a {\em complex (resp. real)
measurable function on $\Sigma$}
we mean a map of abstract $\sigma$-algebras from $\mb(\C)$ (resp. $\mb(\R)$) to
$\Sigma$. Since absolute value is continuous, by Lemma \ref{lb} we
have for a complex measurable function $f$ on $\Sigma$ a non-negative
measurable function $|f|$. If $\Sigma$ is equipped with
a measure, we say that a complex measurable function $f$
is {\em integrable} if 
$$\int |f|<\infty.$$
The $\sigma$-frame of open sets of $\C$ is the coproduct of two copies
of the $\sigma$-frame of open sets of $\R$. Thus, we may write a complex
measurable function on $\Sigma$ as 
$$f=u+iv$$
where $u,v$ are real measurable functions. We then define $$
u^+=max(u,0),\ u^-=min(u,0),
$$
 and when $f$ is integrable,
$$\int fd\mu=\int u^+d\mu - \int u^-d\mu + i\int v^+d\mu -i\int v^-d\mu.$$
Linearity of the integral on complex integrable functions follows from
Theorem \ref{tadd}.

\vspace{3mm}

\begin{theorem}
(Lebesgue dominated convergence theorem - following Rudin \cite{rudin})
Suppose $f_n$ are complex measurable functions on an
abstract $\sigma$-algebra $\Sigma$
such that
$$f=\cform{\lim}{n\r\infty}{}f_n$$
and there exists an integrable non-negative function $g$ on $\Sigma$
such that 
$$|f_n|\leq g.$$
Then $f$ is integrable,
$$\cform{\lim}{n\r\infty}{}\int |f_n-f|d\mu=0$$
and
$$\cform{\lim}{n\r\infty}{}\int f_nd\mu=\int fd\mu.$$
\end{theorem}

\Proof
Apply Fatou's lemma to the function $2g-|f_n-f|$.
\qed

\vspace{3mm}

\section{Abstract radonification}

\label{sradon}

\vspace{3mm}
There is an obvious forgetful functor from 
abstract $\sigma$-algebras to Boolean algebras.
Moreover, the functor obviously preserves limits, and, in
fact, a left adjoint exists by the well known general principle for
(possibly infinitary) unversal algebras. 
We will denote the left adjoint by $CB(?)$.

\vspace{3mm}

\begin{theorem}
\label{trad}
Let $B$ be a Boolean algebra and $\mu_0:B\r[0,1]$ an additive function.
Then there exists a unique measure $\mu:CB(B)\r[0,1]$ such that the
following diagram commutes:
$$
\diagram
B\dto_\eta\rto^{\mu_0}&[0,1]\\
CB(B)\urto_\mu &.
\enddiagram
$$
\end{theorem}

\Proof
By transfinite induction, we shall construct Boolean algebras $B_\alpha$
and additive functions
$$\mu_\alpha:B_\alpha\r[0,1].$$
Put $B_0:=B$. Pick a countable subset 
$$S_\alpha=\{s_{\alpha0},s_{\alpha1},...\}\subset B_\alpha$$
such that 
$$S_\alpha\neq S_\gamma \;\text{for any $\gamma<\alpha$}.$$
Then let 
$$B_{\alpha+1}:=B_{\alpha}\amalg 2^{\{x_\alpha,y_\alpha\}}$$
and define
$$\mu_{\alpha+1}(a\otimes 1):=\mu_{\alpha}(a),$$
\beg{erad+}{\mu_{\alpha+1}(a\otimes x_\alpha):=
\cform{\lim}{n\r\infty}{}\mu_\alpha(a\wedge(s_{\alpha1}\vee...\vee s_{\alpha n})).
}
Additivity of limits clearly implies that this extends uniquely to an additive
function 
$$\mu_{\alpha+1}:B_{\alpha+1}\r[0,1].$$
For a limit ordinal $\alpha$, we define $B_\alpha$ be the union of $B_\beta$
over $\beta<\alpha$, and we let $\mu_\alpha$ be the common extension of 
$\mu_\beta$, $\beta<\alpha$. 

Now let $\gamma$ be a cardinal number $>|B|$ of cofinality $>\omega$
such that 
$$\alpha<\gamma\Rightarrow \alpha^{\aleph_0}<\gamma.$$
Then $\{S_\alpha|\alpha<\gamma\}$ is the set of {\em all} countable subsets
of $B_\gamma$. Let 
$$\widetilde{B}:=B_\gamma/I_{\mu_\gamma}$$
and let $\widetilde{\mu}:\widetilde{B}\r[0,1]$ be the induced additive function.

\vspace{3mm}
\begin{lemma}
\label{lrad}
$\widetilde{B}$ is an
abstract $\sigma$-algebra and $\widetilde{\mu}:\widetilde{B}\r[0,1]$
is a measure.
\end{lemma}

\Proof
First, we claim that for every $\alpha<\gamma$,
\beg{elrad1}{x_\alpha=\bigvee S_\alpha.
}
In effect, 
$$\mu_\gamma(s_{\alpha i}\wedge y_\alpha) =
\mu_{\alpha+1}(s_{\alpha i}\wedge y_\alpha)=
\mu_{\alpha+1}(s_{\alpha i})-\mu_{\alpha+1}(s_{\alpha i}\wedge x_\alpha)=0,$$
and hence $s_{\alpha i}\wedge y_\alpha=0$, so $x_\alpha\geq s_{\alpha i}$.
Now suppose $z\geq s_{\alpha i}$ for all $i=0,1,...$.
Let $u:=x_\alpha\smallsetminus z$. Then $x_\alpha\geq u$ and
$u\wedge s_{\alpha i}=0$ for all $i$. This means that
$$\widetilde{\mu}(u)+\mu_{\alpha}\left(\cform{\bigvee}{i=0}{n}s_{\alpha i}\right)
\leq \mu_{\alpha+1}(x_\alpha)=\cform{\lim}{n\r\infty}{}\mu_\alpha\left(
\cform{\bigvee}{i=0}{\infty}s_{\alpha i}\right)$$
(the last equality is by \rref{erad+}). Hence,
$\widetilde{\mu}(u)=0$ so $u=0$, which proves \rref{elrad1}. Note that
we have also proved
\beg{elrad2}{\cform{\lim}{n\r\infty}{}\widetilde{\mu}\left(\cform{\bigvee}{i=0}{n}
s_{\alpha i}\right)=\widetilde{\mu}\left(\cform{\bigvee}{i=0}{\infty}s_{\alpha i}
\right),
}
which concludes the proof of the Lemma.
\qed

\vspace{3mm}
To continue the proof of the Theorem, we have constructed a diagram
$$\diagram
B\rto^{\mu_0}\dto_\beta & [0,1]\\
\widetilde{B}\urto_{\widetilde{\mu}} &
\enddiagram
$$
where $\widetilde{B}$ is an
abstract $\sigma$-algebra, the inclusion $\beta$ is a morphism
of Boolean algebras, and $\widetilde{\mu}$ is a measure.
By adjunction, we have a unique morphism of 
abstract $\sigma$-algebras $\widetilde{\beta}$
completing the following diagram
$$
\diagram
B\dto\drto^\beta &\\
CB(B) \rdotted|>\tip_{\widetilde{\beta}}&\widetilde{B}.
\enddiagram
$$
Then $\widetilde{\mu}\widetilde{\beta}$ is a measure, proving the existence
statement of the Theorem. To prove uniqueness, note that $\widetilde{B}$
is generated by $B$ as an
abstract $\sigma$-algebra and therefore the measure of each
element of $CB(B)$ can be computed recursively from $\mu_0$ by transfinite 
induction.
\qed

\vspace{3mm}
{\bf Remark:} From our setup, one may perhaps expect that Theorem \ref{trad}
would have a simple proof using the universal property of $CB(?)$. However, at
present we
don't know such a proof. 

\vspace{3mm}

\begin{corollary}
\label{crad1}
Under the assumptions of Theorem \ref{trad}, there exists a
{\em unique} (up to unique isomorphism) morphism of Boolean algebras
$$\beta:B\r \widetilde{B}
$$
where $\widetilde{B}$ is a reduced
abstract $\sigma$-algebra with a measure $\widetilde{\mu}$,
such that $\widetilde{B}$ is generated by $Im(\beta)$, and $\widetilde{\mu}\beta=
\mu_0$.
\end{corollary}

\Proof
For existence, we may take $\widetilde{B}$, $\widetilde{\mu}$ 
as constructed in the proof of Theorem \ref{trad}. For uniqueness, 
if we pick $\widetilde{B}$, $\widetilde{\mu}$ as in the statement
of the Corollary, Theorem \ref{trad} gives a morphism of 
abstract $\sigma$-algebras
\beg{erad4}{\widetilde{\beta}:CB(B)\r \widetilde{B}
}
such that $\widetilde{\mu}\widetilde{\beta}=\mu$.
Additionally, the assumption that $\widetilde{B}$ is generated by
the image of $\beta$ implies that \rref{erad4} is onto. 
Now since $\widetilde{B}$ is reduced, \rref{erad4} factors through
a unique morphism of abstract $\sigma$-algebras
\beg{erad5}{\overline{\beta}:CB(B)/I_\mu\r\widetilde{B}.
}
Then \rref{erad5} is onto because \rref{erad4} is, and
$$\overline{\beta}(x)=0\Rightarrow \mu(x)=\widetilde{\mu}\widetilde{\beta}(x)=0
\Rightarrow x=0,$$
so $\overline{\beta}$ is injective and hence an isomorphism.
\qed

\vspace{3mm}
\noindent
{\bf Remark:}
Note that to construct, say, a non-negative measurable function
on $\widetilde{B}$, all we need is an order-preserving map
$$f:[0,\infty]\r CB(B)$$
(on the left hand side, we identify $t$ with $[0,t]$)
such that for $t_n\searrow t$,
$$\cform{\lim}{n\r\infty}{}\mu(f(t_n))=\mu(f(t)).
$$
Since $\mu$ can be (theoretically) computed from 
$\sigma$-aditivity and $\mu_0$, measurable maps are, in principle,
in abundant supply.

\vspace{3mm}
\noindent
{\bf Example:}
An (irreducible) Gaussian probability space consists of a Hilbert space
$H$, a probability space $(X,\Sigma, \mu)$ and a Hilbert space embedding
\beg{eradi}{\alpha:H\subset L^2(X)} 
such that all the variables in the image of
\rref{eradi} have centered Gaussian law, and they generate $\Sigma$.

\vspace{3mm}
\begin{theorem}
\label{tgauss}
Let $B$ denote the set of cylindrical Borel measurable subsets of $H$, and
let $\mu_0:B\r [0,1]$ be the Gaussian additive function. Then there is a canonical
isomorphism 
$$\widetilde{B}\cong \Sigma/I_\mu.$$
\end{theorem}

\Proof
We define a map of Boolean algebras
$$f:B\r\Sigma$$
as follows; for $a\in H$, $b\in \mb(\R)$ (the Borel $\sigma$-algebra on $\R$),
$$f(\langle a,?\rangle^{-1}(b)):=\alpha^{-1}(b).$$
Clearly, this extends uniquely to a map of Boolean algebras, and
we have a commutative diagram
$$\diagram
B\dto_f\rto^{\mu_0}&[0,1]\\
\Sigma.\urto_\mu &
\enddiagram
$$
Additionally, by definition, $\Sigma$ is 
generated by $f(B)$ as a $\sigma$-algebra
(and hence an
abstract $\sigma$-algebra). 
Thus, our statement follows from Corollary \ref{crad1}.
\qed

\vspace{3mm}
A {\em radonification} of $H$ is a bounded injective dense linear map $\iota$
from $H$ into a Banach space $\overline{H}$, a bounded linear map 
$\overline{?}:H\r \overline{H}$ such that $\overline{a}\iota=\langle a,?\rangle$
and a $\sigma$-algebra $\Sigma$ on $\overline{H}$ with a probability measure $\mu$
such that $\overline{a}$ are measurable functions and for $b\in \mb(\R^n)$,
$$\mu(\cform{\bigcap}{i=1}{n}\overline{a_i}^{-1}(b))=\mu_0(\cform{\bigcap}{i=1}{n}a_{i}^{-1}(b)).$$

\vspace{3mm}
\begin{proposition}
\label{prad}
In this situation, there is a unique isomorphism of abstract $\sigma$-algebras 
$$\widetilde{B}\cong \Sigma/I_\mu$$
which sends $\mu$ to $\tilde{\mu}$.
\end{proposition}

\Proof
$a^{-1}(b)\mapsto \overline{a}^{-1}(b)$ clearly defines an embedding of Boolean
algebras
$$f:B\r\Sigma$$
such that $\mu f=\mu_0$, so again, our statement follows from Corollary \ref{crad1}. \qed

\vspace{3mm}
\noindent
{\bf A non-Gaussian example:} Using Corollary \ref{crad1}, we may obtain an
infinite-dimensional point-free abstract integration theory for any sequence
of absolutely continuous independent random variables $X_1,X_2,$ $\dots$. 
If $H$ is a real Hilbert space with Hilbert basis $e_1,e_2,\dots$, and
if $f_n$ is the probability density of $X_n$, and 
$$\pi_n:H\r \langle e_1,\dots,e_n\rangle$$
is the orthogonal projection, we may define, for a Borel subset $S\subseteq
\langle e_1,\dots,e_n\rangle$,
$$\mu(\pi^{-1}(S))=\int_Sf_1(x_1)\cdot \dots \cdot f_n(x_n)dx_1\dots dx_n.$$
This is a finitely additive function from the set $B$ of cylindrical Borel 
measurable subset of $H$ with values in $[0,1]$, so it extends
to a $\sigma$-additive function $\widetilde{\mu}:\widetilde{B}\r[0,1]$.

Most point-free measures $\widetilde{\mu}$ one obtains in this way are,
however, unnatural in the sense that they are highly dependent on the
coordinate system we choose. It is interesting, in this context, to search
for finite-dimensional
joint distributions which, with respect to onto linear maps, transform
(contravariantly) in families with finitely many parameters in each
dimension. Such non-Gaussian examples are hard to come by.

A beautiful example of this kind is the multivariate skew-normal distribution
introduced by A. Azzalini and A. Dalla Valle \cite{adv}. For the most general
finite-dimensional version of this distribution, consider an
$(n+1)\times (n+1)$ positive-definite real symmetric matrix
$A=(a_{ij})_{0\leq i,j\leq n}$ such that $a_{11}=1$. Let $(X_0,X)^T\in \R\times \R^n$
be a centered normally distributed random variable with covariance matrix $A$.
Then let $Z$ be the random variable give by
$$Z=\left\{
\begin{array}{ll}
X & \text{if $X_0>0$}\\
-X & \text{if $X_0<0$.}
\end{array} \right.$$
The random variable $Z$ is absolutely continuous in $\R^n$ with probability density
function which we will denote by $\phi_A$. This is the skew-normal distribution
(see \cite{adv}, Proposition 6).

For our purposes, it suffices to consider the case when
$(a_{ij})_{1\leq i,j\leq n}$ is the identity matrix, since
one can always reach this case by linear transformation. 
Let, then $\delta_1,\delta_2,\dots \R$ be such that 
\beg{edelta1}{\sum_{n=1}^{\infty}\delta_{n}^{2}<1.
}
Let 
$$A_n=\left(\begin{array}{ccccc}
1&\delta_1&\delta_2&\dots & \delta_n\\
\delta_1&1&0&\dots & 0\\
\delta_2&0&1&\dots & 0\\
\dots &\dots&\dots&\dots&\dots\\
\delta_n&0&0&\dots&1
\end{array}
\right)$$
(\rref{edelta1} assures that $A_n$ is positive-definite). Define an additive
function with values in $[0,1]$ on the set $B$ of cylindrical Borel-measurable
subsets of the Hilbert space $H$ by setting, for a Borel-measurable subset
$S\subseteq \langle e_1,\dots,e_n\rangle$, 
$$\mu(\pi_{n}^{-1}S)=\int_S \phi_{A_n}(x_1,\dots,x_n)dx_1\dots dx_n.$$
Applying Corollary \ref{crad1}, we obtain a reduced abstract $\sigma$-algebre
$\widetilde{B}$ with a measure $\widetilde{\mu}$. This is the infinite-dimensional
skew-normal distribution. As far as we know, a pointed version of this
distribution has not been considered.

\end{document}